\documentclass[reqno]{amsart}
\usepackage{hyperref}

\AtBeginDocument{{\noindent\small
Original Russian Text \copyright \, T.N.~Harutyunyan, A.A.~Pahlevanyan, Yu.A.~Ashrafyan, 2013 \\
published in \emph{Proceedings of Artsakh State University}, (2013), No. 1(27), pp. 12--20.}
\vspace{9mm}}

\begin{document}
\title[On the "movement" of the zeros of eigenfunctions of the S-L problem]{On the "movement" of the zeros of eigenfunctions of the Sturm-Liouville problem}

\author{Tigran Harutyunyan, Avetik Pahlevanyan, Yuri Ashrafyan}

\address{Tigran Harutyunyan \newline
Faculty of Mathematics and Mechanics, Yerevan State University, 1 Alex Manoogian, 0025, Yerevan, Armenia}
\email{hartigr@yahoo.co.uk}
\address{Avetik Pahlevanyan \newline
Faculty of Mathematics and Mechanics, Yerevan State University, 1 Alex Manoogian, 0025, Yerevan, Armenia}
\email{avetikpahlevanyan@gmail.com}
\address{Yuri Ashrafyan \newline
Faculty of Mathematics and Mechanics, Yerevan State University, 1 Alex Manoogian, 0025, Yerevan, Armenia}
\email{yuriashrafyan@ysu.am}


\begin{abstract}
We study the dependence of the zeros of eigenfunctions of Sturm-Liouville problem on the parameters that define the boundary conditions. As a corollary, we obtain Sturm oscillation theorem, which states that the $n$-th eigenfunction has $n$ zeros.
\end{abstract}

\maketitle
\newtheorem{theorem}{Theorem}
\newtheorem{lemma}[theorem]{Lemma}
\newtheorem{definition}[theorem]{Definition}
\newtheorem{remark}[theorem]{Remark}

Let us consider Sturm-Liouville boundary value problem $L\left(q, \alpha, \beta \right):$
\begin{equation}\label{eq1}
ly \equiv - y'' + q\left( x \right)y = \mu y, \; 0<x<\pi, \; \mu \in \mathbb{C},
\end{equation}
\begin{equation}\label{eq2}
y\left( 0 \right)\cos \alpha  + y'\left( 0 \right)\sin \alpha  = 0, \; \alpha  \in \left( {0,\pi } \right],
\end{equation}
\begin{equation}\label{eq3}
y\left( \pi \right)\cos \beta  + y'\left( \pi \right)\sin \beta  = 0, \; \beta  \in \left[ {0,\pi } \right),
\end{equation}
where $q \in L_{\mathbb{R}}^1\left[ {0,\pi } \right],$ i.e. $q$ is a real-valued, summable function on $\left[0, \pi\right]$.

By $L\left(q,\alpha ,\beta \right)$ we also denote the self-adjoint operator corresponding to the problem \eqref{eq1}--\eqref{eq3}.

It is well-known, that the problem $L\left(q, \alpha, \beta \right)$ has a countable set of simple, real eigenvalues (see, e.g. \cite{Levitan_Sargsyan:1988,Marchenko:1977,Freiling_Yurko:2001,Harutyunyan:2008}), which we denote by ${\mu }_n\left(q,\alpha ,\beta \right),$ $n=0,1,\dots,$ (emphasizing the dependence on $q,$ $\alpha$ and $\beta$) and enumerate in increasing order:
\begin{equation}\label{eq4}
\mu_0 \left(q,\alpha,\beta\right) < \mu_1 \left(q,\alpha,\beta\right) < \dots < \mu_n \left(q,\alpha,\beta\right) < \dots \; .
\end{equation}
In the papers \cite{Harutyunyan:2008,Harutyunyan_Navasardyan:2000} it was introduced the concept of the eigenvalues function (EVF) of a family of Sturm-Liouville operators $\left\{L\left(q, \alpha, \beta\right); \alpha \in \left(0,\pi\right], \beta\left[0,\pi\right)\right\}.$ For fixed $q$ it is a function of two variables $\gamma$ and $\delta$ $\left(\gamma=\alpha+\pi n \in \left(0,\infty\right), \delta=\beta-\pi m \in \left(-\infty, \pi\right)\right)$ determined through the eigenvalues $\mu_n \left(q,\alpha,\beta\right),$ $n=0,1,2,\dots,$ by the following formula:
\begin{equation}\label{eq5}
\mu \left(\gamma,\delta\right)=\mu\left(\alpha+\pi n, \beta-\pi m \right):= \mu_{n+m} \left(q,\alpha,\beta\right), \; n,m=0,1,2,\dots.
\end{equation}
It is proved that this function is analytic with respect to $\gamma$ and $\delta.$ It is strictly increasing by $\gamma$ and strictly decreasing by $\delta.$

It is also known (see, e.g. \cite{Levitan_Sargsyan:1988}) that every nontrivial solution $y\left(x,\mu\right)$ of the equation \eqref{eq1} may have only simple zeros (if $y\left(x_0,\mu\right)=0,$ then $y'\left(x_0,\mu\right)\neq 0$), and (see, e.g. \cite{Ghazaryan:2002}) every solution $y\left(x,\mu\right)$ is a continuously differentiable function with respect to the totality of variables $x$ and $\mu.$ Therefore, by applying the implicit function theorem (see, e.g. \cite[p. 452]{Fikhtengolts:1966}), we get that the zeros of the solution $y\left(x,\mu\right)$ are the continuously differentiable functions with respect to $\mu.$ Since the function $x=x\left(\mu\right),$ such that the identity $y\left(x\left(\mu\right),\mu\right) \equiv 0$ is true for all $\mu$ from some interval $\left(a,b\right),$ is called the solution of the equation $y\left(x,\mu\right)=0,$ then differentiating the last identity with respect to $\mu$ we obtain:
\begin{equation}\label{eq6}
\cfrac{dy\left(x\left(\mu\right), \mu\right)}{d\mu}=\cfrac{\partial y\left(x\left(\mu\right), \mu\right)}{\partial x} \, \cfrac{dx\left(\mu\right)}{d\mu}+\cfrac{\partial y\left(x\left(\mu\right), \mu\right)}{\partial \mu} \equiv 0, \; \mu \in \left(a,b\right).
\end{equation}
Let us denote $\cfrac{\partial y\left(x,\mu\right)}{\partial \mu}:=\dot{y}\left(x,\mu\right)$ and write the identity \eqref{eq6} in the following form:
\begin{equation}\label{eq7}
\cfrac{dx\left(\mu\right)}{d\mu}=\dot{x}\left(\mu\right)=-\cfrac{\dot{y}\left(x\left(\mu\right), \mu\right)}{y'\left(x\left(\mu\right), \mu\right)}, \; \mu \in \left(a,b\right).
\end{equation}
On the other hand, let us write down the fact that $y\left(x,\mu\right)$ is the solution of the equation \eqref{eq1}, i.e.
\begin{equation}\label{eq8}
-y''\left(x, \mu \right)+q\left(x\right)y\left(x,\mu \right)\equiv \mu y\left(x,\mu \right), \; 0<x<\pi, \; \mu \in \mathbb{C},
\end{equation}
and differentiating this identity with respect to $\mu$ we receive:
\begin{equation}\label{eq9}
-\dot{y}'' \left(x, \mu \right)+q\left(x\right)\dot{y}\left(x, \mu \right) \equiv y\left(x, \mu \right)+\mu \dot{y}\left(x,\mu \right).
\end{equation}
Multiplying \eqref{eq8} by $\dot{y},$ \eqref{eq9} by $y$ and subtracting from the second obtained identity the first one, we get
\begin{equation*}
y''\left(x, \mu \right)\dot{y}\left(x, \mu \right)-\dot{y}''\left(x, \mu \right)y\left(x, \mu \right) \equiv y^{2} \left(x, \mu \right), \; 0<x<\pi, \; \mu \in {\mathbb C},
\end{equation*}
i.e.
\begin{equation}\label{eq10}
\cfrac{d}{dx} \left[y'\left(x, \mu \right)\dot{y}\left(x, \mu \right)-\dot{y}'\left(x, \mu \right)y\left(x, \mu \right)\right] \equiv y^{2} \left(x, \mu \right).
\end{equation}
If we integrate this identity with respect to $x$ from $0$ to $a$ $\left(0 \leq a \leq \pi \right),$ then
\begin{multline}\label{eq11}
y'\left(a, \mu \right)\dot{y}\left(a, \mu \right)-\dot{y}'\left(a, \mu \right)y\left(a, \mu \right)-y'\left(0, \mu \right)\dot{y}\left(0, \mu \right)+\\
+\dot{y}'\left(0, \mu \right)y\left(0, \mu \right)=\int_{0}^{a} y^{2} \left(x, \mu \right)dx,
\end{multline}
and if we integrate with respect to $x$ from $a$ to $\pi,$ then we get
\begin{multline}\label{eq12}
y'\left(\pi, \mu \right)\dot{y}\left(\pi, \mu \right)-\dot{y}'\left(\pi, \mu \right)y\left(\pi, \mu \right)-y'\left(a, \mu \right)\dot{y}\left(a, \mu \right)+\\
+\dot{y}'\left(a, \mu \right)y\left(a, \mu \right)=\int_{a}^{\pi} y^{2} \left(x, \mu \right)dx.
\end{multline}
Now, as $y\left(x,\mu\right)$ let us take $y=\varphi\left(x, \mu, \alpha, q\right)-$ the solution of the equation \eqref{eq1}, satisfying the following initial conditions:
\begin{equation}\label{eq13}
\varphi \left(0, \mu, \alpha, q\right)=\sin \alpha, \;\;\; \varphi'\left(0, \mu, \alpha, q\right)=-\cos \alpha.
\end{equation}
It is easy to see that eigenfunctions of the problem $L\left(q,\alpha,\beta\right)$ are obtained from the solution $\varphi\left(x, \mu, \alpha, q\right)$ at $\mu=\mu_n\left(q,\alpha,\beta\right)$ (here we use \eqref{eq5}), i.e.
\begin{multline}\label{eq14}
\varphi_{n} \left(x, q, \alpha, \beta \right):=\varphi_{n} \left(x\right)=\varphi \left(x, \mu_{n} \left(q, \alpha, \beta \right), \alpha, q\right)=\\
=\varphi \left(x, \mu \left(\alpha+\pi n, \beta \right), \alpha, q\right)=\varphi \left(x, \mu \left(\alpha, \beta-\pi n\right), \alpha, q\right)=\\
=\left. \varphi \left(x, \mu \left(\alpha, \delta\right), \alpha, q\right)\right|_{\delta=\beta-\pi n} =\left. \varphi \left(x, \mu \left(\alpha, \delta \right)\right)\right|_{\delta =\beta-\pi n} =\\
=\varphi \left(x, \mu \left(\alpha, \beta -\pi n\right)\right).
\end{multline}
Let $0 \leq x_{n}^{0} < x_{n}^{1} < \dots < x_{n}^{m} \leq \pi$ be the zeros of the eigenfunction $\varphi_{n} \left(x, q, \alpha, \beta \right)=\varphi \left(x, \mu \left(\alpha, \beta-\pi n\right)\right),$ i.e. $\varphi_{n} \left(x_{n}^{k}, q, \alpha, \beta \right)=\varphi \left(x_{n}^{k}, \mu \left(\alpha, \beta-\pi n \right) \right)=0,$ $k=0,1,\dots,m.$

Let $q, \alpha, n$ be fixed. We consider the following questions:
\begin{itemize}
\item[a)]
How are the zeros $x_{n}^{k}=x_{n}^{k} \left(\beta \right),$ $k=0,1,\dots,m$ changing, when $\beta$ is changing on $\left[0,\pi\right)?$

\item[b)]
How many zeros the $n$-th eigenfunction $\varphi _{n} \left(x,q,\alpha ,\beta \right)$ has, i.e. what is equal $m?$
\end{itemize}
By taking $y=\varphi_n\left(x\right)$ and $a=x_{n}^{k}$ $\left(k=0,1,\dots,m\right)$ in \eqref{eq11}, we receive
\begin{multline}\label{eq15}
\varphi'_{n}\left(x_{n}^{k} \right)\dot{\varphi}_{n}\left(x_{n}^{k} \right)-\dot{\varphi}'_{n}\left(x_{n}^{k} \right)\varphi_{n}\left(x_{n}^{k} \right)-\varphi'_{n}\left(0 \right)\dot{\varphi}_{n}\left(0 \right)+\\
+\dot{\varphi}'_{n}\left(0 \right)\varphi_{n}\left(0 \right)=\int_{0}^{x_{n}^{k}} \varphi_{n}^{2} \left(x\right)dx.
\end{multline}
Since the initial conditions \eqref{eq13} should be held for all $\mu \in \mathbb{C},$ then $\dot{\varphi}_{n} \left(0\right)=0$ and $\dot{\varphi}'_{n} \left(0\right)=0.$ Also, taking into account that $\varphi_{n} \left(x_{n}^{k} \right)=0,$ from \eqref{eq15} we get
\begin{equation}\label{eq16}
\varphi'_{n} \left(x_{n}^{k} \right)\dot{\varphi}_{n} \left(x_{n}^{k} \right)=\int_{0}^{x_{n}^{k}} \varphi_{n}^{2} \left(x\right)dx.
\end{equation}
Since the zeros of the solutions are simple, then $\varphi'_{n} \left(x_{n}^{k} \right) \neq 0,$ and therefore, from \eqref{eq16} the following equality implies :
\begin{equation}\label{eq17}
\cfrac{\dot{\varphi}_{n} \left(x_{n}^{k} \right)}{\varphi'_{n} \left(x_{n}^{k} \right)} =\cfrac{1}{\left(\varphi'_{n} \left(x_{n}^{k} \right)\right)^2} \int_{0}^{x_{n}^{k}} \varphi_{n}^{2}\left(x\right)dx.
\end{equation}
Now, from \eqref{eq7}, we obtain
\begin{equation}\label{eq18}
\dot{x}_{n}^{k} \left(\mu_{n} \right)=\left. \cfrac{dx_{n}^{k} \left(\mu \right)}{d\mu} \right|_{\mu=\mu_{n}} =-\cfrac{\dot{\varphi}_{n} \left(x_{n}^{k} \right)}{\varphi'_{n} \left(x_{n}^{k} \right)} =-\cfrac{1}{\left(\varphi'_{n} \left(x_{n}^{k}\right)\right)^2} \int_{0}^{x_{n}^{k}} \varphi_{n}^{2}\left(x\right)dx,
\end{equation}
i.e. zeros $x_{n}^{k} \left(\mu_{n} \right),$ $k=0,1,\dots,m$ of the eigenfunction $\varphi_{n} \left(x\right)$
are decreasing if the eigenvalue $\mu_{n} \left(q, \alpha, \beta \right)$ is increasing, which in its turn means that
\begin{equation}\label{eq19}
\dot{x}_{n}^{k} \left(\mu_{n} \left(q, \alpha, \beta \right)\right) \leq 0.
\end{equation}

Let us note, that the equality $\dot{x}_{n}^{k} \left(\mu_{n} \right)=0$ may occur only at $x_{n}^{k}=0,$ i.e., when $x=0$ is a zero of the eigenfunction $\varphi_{n} \left(x\right),$ and it is so at $\alpha =\pi$  $\left(\gamma=\pi l, \; l=1,2,\dots\right).$

Meanwhile, in the inequality $\dot{x}_{n}^{k} \left(\mu_{n} \right)=\dot{x}_{n}^{k} \left(\mu_{n} \left(q, \alpha, \beta \right)\right) \leq 0,$ the variable $\mu_{n}$ can be changed depending on $q,$ $\alpha$ and $\beta.$ More precisely, if for some change in these three variables $\mu_{n} \left(q, \alpha, \beta \right)$ increases, then the zeros of the eigenfunction $\varphi_{n} \left(x\right)$ are moving to the left, and if $\mu_{n} \left(q, \alpha, \beta \right)$ decreases, then the zeros of the eigenfunction $\varphi_{n} \left(x\right)$ are moving to the right.

In the work \cite{Poschel_Trubowitz:1987}, it was proved for $q \in L_{\mathbb{R}}^{2}\left[0,\pi\right]$ that the number of zeros of the $n$-th eigenfunctions of the problems $L\left(q, \pi, 0 \right)$ and $L\left(0, \pi, 0 \right)$ are equal. But it is easy to see that the same proof remains true for $q \in L_{\mathbb{R}}^{1}\left[0,\pi\right].$

It is easy to calculate that the eigenvalues of the problem $L\left(0,\pi,0\right)$ are $\mu_n\left(0,\pi,0\right)=\left(n+1\right)^2$ and the eigenfunctions are
\begin{multline}\label{eq20}
\varphi_{n} \left(x\right)=\varphi \left(x, \mu_{n} \left(0, \pi, 0 \right), \pi, 0 \right)=\\
=\varphi \left( x, \left(n+1\right)^{2}, \pi, 0 \right)=\cfrac{\sin \left(n+1\right)x}{n+1}, \; n=0,1,2,\dots.
\end{multline}

The zeros of this eigenfunction are $x_{n}^{k}=\cfrac{\pi k}{n+1}, k=0,1,\dots,n+1,$ i.e. the $n$-th eigenfunction of the problem $L\left(0,\pi,0\right)$ has $n+2$ zeros in $\left[0,\pi\right],$ two of which are $0$ and $\pi,$ and $n$ zeros in $\left(0,\pi\right).$

Thereby, the $n$-th eigenfunction $\varphi \left(x, q, \pi, 0 \right)$ of the problem $L\left(q, \pi, 0\right)$ has two zeros at the endpoints of $\left[0,\pi\right],$ i.e. $x_{n}^{0} \left(q, \pi, 0\right)=0,$ $x_{n}^{n+1} \left(q, \pi, 0\right)=\pi$ and also $n$ zeros in the interval $\left(0, \pi\right).$

Along with increasing $\beta$ from $0$ to $\pi$ the eigenvalue $\mu_n \left(q,\pi,0\right)$ is continuously (with respect to $\beta$) decreasing from $\mu_n \left(q,\pi,0\right)$ to $\mu_n \left(q,\pi,\pi\right)=\mu\left(\pi,\pi-\pi n\right)=$ \\
$=\mu\left(\pi,0-\left(n-1\right)\pi\right)=\mu_{n-1} \left(q,\pi,0\right)$ (see \cite{Ghazaryan:2002}) and, according to \eqref{eq18}, the zeros of the function $\varphi_n\left(x,q,\pi,\beta\right)$ are increasing, i.e. are moving to the right (all but leftmost zero $x_n^0=0$). In particular, the rightmost zero $x_n^{n+1}=\pi,$ by moving to the right, leaves the segment $\left[0,\pi\right]$ and in $\left[0,\pi\right]$ remains $n+1$ zeros (one $x_n^0=0$ and $n$ zeros in the interval $\left(0,\pi\right)$). And the previous zero reaches $\pi$ when the relation $\varphi_n\left(\pi\right)=\varphi\left(\pi,\mu_n \left(q,\pi,\beta\right),\pi,q\right)=c_n \psi_n \left(\pi\right) \sin \beta=0$ again occurs (see below \eqref{eq21} and \eqref{eq24}), and this is possible only when $\beta$ reaches $\pi$ (and $\mu_n \left(q,\pi,\beta\right),$ by decreasing, reaches $\mu_n \left(q,\pi,\pi\right)=\mu_{n-1} \left(q,\pi,0\right)$). Then the eigenfunction $\varphi_n \left(x\right)=\varphi\left(x,\mu_n \left(q,\pi,\beta\right),\beta,q\right)$ will smoothly transform to the eigenfunction $\varphi\left(x,\mu_n \left(q,\pi,\pi \right),\pi,q\right)=\varphi\left(x,\mu_{n-1} \left(q,\pi,0 \right),\pi,q\right)$ which has $n+1$ zeros in $\left[0,\pi \right],$ two out of which are the endpoints $0$ and $\pi,$ and $n-1$ zeros are in $\left(0,\pi \right).$ Thus, the oscillation theorem is proved for all $L\left(q,\pi,\beta\right),$ $\beta \in \left[0,\pi\right).$

Now, as $y\left(x,\mu\right)$ let us take $\psi\left(x, \mu, \beta, q\right)-$ the solution of the equation \eqref{eq1} satisfying the following initial conditions:
\begin{equation}\label{eq21}
\psi \left(\pi, \mu, \beta, q\right)=\sin \beta, \;\;\; \psi'\left(\pi, \mu, \beta, q\right)=-\cos \beta.
\end{equation}
It is easy to see that the eigenvalues $\mu_n=\mu_n \left(q,\alpha,\beta\right),$ $n=0,1,\dots,$ of the problem $L\left(q,\alpha,\beta\right)$ are the zeros of the entire function
\begin{equation}\label{eq22}
\Psi \left(\mu \right)=\Psi \left(\mu, \alpha, \beta, q\right)=\psi \left(0, \mu, \beta, q\right)\cos \alpha +\psi'\left(0, \mu, \beta, q\right)\sin \alpha,
\end{equation}
and the eigenfunctions, corresponding to these eigenvalues, are obtained by the formula
\begin{equation}\label{eq23}
\psi_{n} \left(x\right)=\psi \left(x, \mu_{n} \left(q, \alpha, \beta \right), \beta, q\right), \; n=0,1,\dots.
\end{equation}
Since all eigenvalues $\mu_n$ are simple, then eigenfunctions $\varphi_n \left(x\right)$ and $\psi_n \left(x\right)$ corresponding to the same eigenvalue $\mu_n$ are linearly dependent, i.e. there exist the constants $c_n=c_n \left(q,\alpha,\beta\right),$ $n=0,1,\dots,$ such that
\begin{equation}\label{eq24}
\varphi_{n} \left(x\right)=c_{n} \psi_{n} \left(x\right), \; n=0,1,\dots.
\end{equation}
This implies that $\varphi_{n} \left(x,q,\alpha,\beta \right)$ and $\psi_{n} \left(x,q,\alpha,\beta \right)$ have equal number of zeros.

Let $0 \leq x_n^m < x_n^{m-1} < \dots <x_n^0 \leq \pi$ be the zeros of the eigenfunction $\psi_{n} \left(x,q,\alpha,\beta \right)$ i.e. $\psi_{n} \left(x_n^k, q, \alpha, \beta \right)=0,$ $k=0,1,\dots,m.$

By taking in the identity \eqref{eq12} $y=\psi_{n} \left(x\right)=\psi \left(x, \mu_n \left(q, \alpha, \beta \right), \beta, q \right),$ we receive
\begin{multline}\label{eq25}
\psi'_{n} \left(\pi \right)\dot{\psi}_{n} \left(\pi \right)-\dot{\psi}'_{n} \left(\pi \right)\psi_{n} \left(\pi \right)-\psi'_{n} \left(x_{n}^{k} \right)\dot{\psi}_{n} \left(x_{n}^{k} \right)+\\
+\dot{\psi}'_{n} \left(x_{n}^{k} \right)\psi_{n} \left(x_{n}^{k} \right)=\int_{x_{n}^{k}}^{\pi}\psi_{n}^{2} \left(x\right)dx.
\end{multline}
From \eqref{eq21} we get that $\dot{\psi}_{n} \left(\pi\right)=0$ and $\dot{\psi}'_{n} \left(\pi\right)=0.$ And since $\psi_{n} \left(x_n^k\right)=0,$ then the equality \eqref{eq25} gets the form
\begin{equation}\label{eq26}
-\psi'_{n} \left(x_{n}^{k} \right)\dot{\psi}_{n} \left(x_{n}^{k} \right)=\int_{x_{n}^{k}}^{\pi} \psi_{n}^{2} \left(x\right)dx.
\end{equation}
As all zeros $x_n^k$ are simple, i.e. $\psi'_{n} \left(x_n^k\right) \neq 0,$ dividing both sides of the last equality by $\left(\psi'_{n} \left(x_n^k\right)\right)^2,$ we obtain
\begin{equation}\label{eq27}
\cfrac{\dot{\psi}_{n} \left(x_{n}^{k} \right)}{\psi'_{n} \left(x_{n}^{k} \right)} =-\cfrac{1}{\left(\psi'_{n} \left(x_{n}^{k} \right)\right)^2} \int_{x_{n}^{k}}^{\pi} \psi_{n}^{2} \left(x\right)dx.
\end{equation}
Now from \eqref{eq7}, by taking $y=\psi_n \left(x\right),$ $x=x_n^k,$ we have
\begin{equation}\label{eq28}
\dot{x}_{n}^{k} \left(\mu_{n} \right)=\left. \cfrac{dx_{n}^{k} \left(\mu \right)}{d\mu} \right|_{\mu =\mu_{n}} =-\cfrac{\dot{\psi}_{n} \left(x_{n}^{k} \right)}{\psi'_{n} \left(x_{n}^{k} \right)} =\cfrac{1}{\left(\psi'_{n} \left(x_{n}^{k} \right)\right)^2} \int_{x_{n}^{k}}^{\pi} \psi_{n}^{2} \left(x\right)dx \geq 0,
\end{equation}
i.e. the zeros $x_n^k \left(\mu_n\right),$ $k=0,1,\dots,m$ of the eigenfunction $\psi_n \left(x\right)$ are increasing, if the eigenvalue $\mu_n$ is increasing. Note that the equality $\dot{x}_n^k \left(\mu_n\right)=0$ is possible only when $x_n^0 \left(\mu_n\right)=\pi,$ but it holds when $\beta=0$ $\left(\delta=-\pi l, \; l=0,1,2,\dots\right).$

While studying the dependence of the zeros of eigenfunctions on $\alpha,$ it is convenient to use formula \eqref{eq28}, because the eigenfunctions $\psi_n \left(x\right)$ have fixed values $\psi_n \left(\pi, \mu, \beta\right)=\sin \beta$ and $\psi'_n \left(\pi, \mu, \beta\right)=-\cos \beta$ for all $\mu \in \mathbb{C},$ i.e. all $\psi_n \left(x\right)$ satisfy the initial conditions $\psi_n \left(\pi\right)=\sin \beta$ and $\psi'_n \left(\pi\right)=-\cos \beta,$ which means that from the endpoint $\pi$ of the segment $\left[0,\pi\right]$ (when changing $\alpha$) new zeros can neither enter nor leave (neither appear nor disappear). Thus, with increasing $\alpha,$ the eigenvalues $\mu_n \left(q, \alpha, \beta\right)$ (with fixed $q$ and $\beta$) are increasing and according to \eqref{eq28} (i.e. $\dot{x}_n^k \left(\mu_n\right) \geq 0$) the zeros of the eigenfunction $\psi_n \left(x\right)$ are moving to the right (i.e. are increasing). Wherein, the values $\psi \left(\pi\right)=\sin \beta$ and $\psi' \left(\pi\right)=-\cos \beta$ are not changed, the number of the zeros increase, and these zeros can neither ``collide'' nor ``split'' as they are simple. That's why new zeros can appear only entering through $0-$ left endpoint of the segment $\left[0,\pi\right]$ and moving to the right (and ``condensing'' respectively). And new zeros are entered through $0-$ left endpoint of the segment $\left[0,\pi\right]$ only when $\psi_n \left(0\right)=0,$ but since $\psi_n \left(0\right)=c_n \varphi_n \left(0\right)=c_n \sin \alpha$ $\left(c_n \neq 0\right),$ then the equality $\psi_n \left(0\right)=0$ is possible only when $\sin \alpha=0,$ i.e. in our notations only when $\alpha=\pi$ (and because $\mu_n \left(q, 0, \beta\right)=\mu \left(0+\pi n, \beta\right)=\mu \left(\pi+\left(n-1\right)\pi, \beta\right)=\mu_{n-1} \left(q, \pi, \beta\right)$) or when $\alpha=0.$

Hence, when $\alpha=\pi,$ the eigenfunction $\psi_n \left(x, q, \pi, \beta\right)$ as well as $\varphi_n \left(x, q, \pi, \beta\right)$ has $n$ zeros in $\left(0,\pi\right)$ and one zero at $x=0$ left endpoint. Wherein, $\psi_n \left(x, q, \pi, \beta\right)=\psi \left(x, \mu_n \left(q, \pi, \beta\right), \beta, q\right)=\psi \left(x, \mu_{n+1} \left(q, 0, \beta\right), \beta, q\right)=\psi_{n+1} \left(x, q, 0, \beta\right).$ With increasing $\alpha$ from $0$ to $\pi$ the eigenvalue $\mu_{n+1} \left(q, 0, \beta\right)$ is increasing (continuously with respect to $\alpha$) to $\mu_{n+1} \left(q, \pi, \beta \right),$ and the leftmost zero $x=0$ by moving to the right, appears in $\left(0,\pi\right),$ i.e. there are $n+1$ zeros of the eigenfunction $\psi_n \left(x, q, \alpha, \beta\right)$ in $\left(0,\pi\right)$ (and another fixed zero $x=\pi,$ if $\beta=0$). A new zero will appear at $x=0$ left endpoint, when $\alpha$ will reach $\alpha=\pi.$

Thus, we obtain the following oscillation theorem:
\begin{theorem}\label{thm1}
The eigenfunctions of the problem $L\left(q, \alpha, \beta\right)$ corresponding to the $n$-th eigenvalue $\mu_n \left(q, \alpha, \beta\right),$ $n=0,1,2,\dots,$ have exactly $n$ zeros in $\left(0,\pi\right).$ All these zeros are simple. If $\alpha=\pi$ and $\beta=0,$ then the $n$-th eigenfunction has $n+2$ zeros in $\left[0,\pi\right],$ and if either $\alpha=\pi,$ $\beta \in \left(0,\pi\right)$ or $\beta=0,$ $\alpha \in \left(0,\pi\right),$ then the $n$-th eigenfunction has $n+1$ zeros in $\left[0,\pi\right].$
\end{theorem}
Oscillation properties of the solutions of the problem $L\left(q, \alpha, \beta\right)$ (Sturm theory), the studies of which initiated by Sturm in \cite{Sturm1:1836,Sturm2:1836}, outlined in monographic literature (see, e.g. \cite{Levitan_Sargsyan:1988,Sansone:1953,Coddington_Levinson:1955}) for continuous $q.$ In the recent years, the study was mostly focused on the cases when $q$ is bounded or $q \in L_{\mathbb{R}}^{2} \left[0,\pi\right]$ (see, e.g. \cite{Poschel_Trubowitz:1987,Simon:2005}), but in many studies (see \cite{Hinton:2005} and references therein) implicitly assumes that Sturm's oscillation theorem (that the $n$-th eigenfunction has $n$ zeros) is also true for $q \in L_{\mathbb{R}}^{1} \left[0,\pi\right],$ although the rigorous proof is not available in the literature.

Our oscillation theorem is true for all $q \in L_{\mathbb{R}}^{1} \left[0,\pi\right].$

\end{document}